\documentclass  [12pt]{amsproc}

\newtheorem{theorem}{Theorem}[section]
\newtheorem{lemma}[theorem]{Lemma}

\newtheorem{definition}[theorem]{Definition}

\newtheorem{corollary}[theorem]{Corollary}

\theoremstyle{remark}

\numberwithin{equation}{section}

\begin{document}
\title {Structured Parseval Frames in Hilbert $C^*$-modules}

\author{Wu Jing}
\address{Department of Mathematics, University of central Florida,  Orlando, FL 32816}
\email{mathjing@gmail.com}

\author{Deguang Han}
\address{Department of Mathematics, University of central Florida,  Orlando, FL 32816}
\email{dhan@pegasus.cc.ucf.edu}

\author{Ram N. Mohapatra}
\address{Department of Mathematics, University of central Florida,  Orlando, FL 32816}
\email{ramm@mail.ucf.edu}

\subjclass{Primary 46L99; Secondary 42C15, 46H25}

\date{ }
\keywords{Frames, Parseval frame vectors, mutil-frame
approximations, unitary groups, Hilbert $C^*$-modules}
\begin{abstract}
We investigate the structured frames for Hilbert $C^{*}$-modules.
In the case that the underlying $C^{*}$-algebra is a commutative
$W^*$-algebra, we prove that the set of the Parseval frame
generators for a unitary operator group  can be parameterized by
the set of all the unitary operators in the double commutant of
the group. Similar result holds for the set of all the general
frame generators where the unitary operators are replaced by
invertible and adjointable operators. Consequently, the set of all
the Parseval frame generators is path-connected. We also obtain
the existence and uniqueness results for the best Parseval
multi-frame approximations for multi-frame generators of unitary
operator groups on Hilbert $C^*$-modules when the underlying
$C^{*}$-algebra is commutative.
\end{abstract}
\maketitle

\section {Introduction}

Frames (modular frames) for Hilbert $C^*$-modules were introduced
by Frank and Larson and some basic properties were also
investigated in a series of their papers \cite{frla207, frla325,
frla273}.  It should be remarked that  although (at the first
glance) some of the definitions and result statements of modular
frames may appear look like similar to their Hilbert space frame
counterparts, these are not simple generalizations of the Hilbert
space frames due to the complexity of the Hilbert $C^{*}$-module
structures and to the fact that many useful techniques in Hilbert
spaces are either not available or not known in Hilbert
$C^{*}$-modules. For example, it is well-known that that every
Hilbert space has  an orthonormal basis which can be simply
obtained by applying the Gram-Schmidt orthonormalization process
to a linearly independent generating subset of the Hilbert space.
However, it is well-known that not every Hilbert $C^{*}$-module
has an ``orthonormal basis". This makes frames particulary
relevant to Hilbert $C^{*}$-modules. Remarkably every countably
generated Hilbert $C^{*}$-module admits a (countable) frame.  It
requires a very deep Hilbert C*-module result (Kasparov's
Stabilization Theorem) to prove this fact which is a trivial fact
in Hilbert space setting (cf \cite{frla273, rath1557}. In fact, it
is still an interesting question whether there exists a
alternative proof for this fact without using the Kasparov's
Stabilization Theorem. Another still open problem is whether every
uncountably generated Hilbert $C^{*}$-module admits a (uncountably
indexed) Parseval frame (again, a trivial fact for Hilbert
spaces). Equivalently, for every Hilbert C*-module over a unital
C*-algebra $\mathcal A$, does there exist an isometric embedding
into a standard Hilbert C*-module $l^{2}({\mathcal A},I)$ as an
orthogonal direct summand for some index set $I$? (see
\cite{frla273}).

In recent years, there have been growing evidence indicating that
 modular frames  are also closely related to some other areas
of research such as the area of wavelet frame constructions (cf.
\cite{Packer04, packer-rieffel03, packer-rieffel04, Wood04}).
Considering the fact that the theory and applications of
structured frames (such as Gabor frame, wavelet frames and frames
induced from group unitary representations) for Hilbert spaces
have been the main focus of the Hilbert space frame theory, we
believe that structure modular frames may well be suitable for
some applications either in theoretical or applied nature. The
purpose of this paper is to initiate the study of structured
modular frames. It is reasonable that we should first take a close
look at
 those existed results for structured Hilbert space frames and make an effort to
 check whether they are still valid for structured modular frames.
 The two results (Theorems \ref{the-3.4} and  \ref{mul-app}) presented in this paper
 are generalizations of the corresponding Hilbert space frame
 results obtained in \cite{hala2000} and  \cite{ha3329}.
 Theorem \ref{the-3.4} states that all the Parseval frame generators for a
 unitary group can be parameterized in terms of the unitary
 elements in the double commutant of the group under the
 commutativity condition on the underlying $C^*$-algebras. This is
 slightly different from the Hilbert space setting since the  the double commutant theorem
 for von Neumann algebras is not always available for the Hilbert $C^{*}$-module
 setting. For the similar reason, the ``finiteness " of the involved
 ``commutant" algebras need to be verified, and each step needs
 to be carefully checked to make sure it is valid in the
 $C^{*}$-algebra context. Theorem \ref{mul-app} deals with the best
 approximations of modular frame generators by Parseval frame
 generators. The difficulty arises when comes to compare two
 positive elements in the underlying $C^{*}$-algebras which is not
 an issue in the scalar case. We are not able to prove these results
 when the underlying $C^{*}$-algebras are non-commutative.

\section{Preliminaries}

This section contains some basic definitions about Hilbert
$C^{*}$-modules and some simple properties for Hilbert
$C^{*}$-module frames that will be needed in the next two
sections. Let $\mathcal A$ be a $C^*$-algebra and ${\mathcal H}$
be a (left) $\mathcal A$-module. Suppose that the linear
structures given on $\mathcal A$ and ${\mathcal H}$ are
compatible, i.e. $\lambda (ax)=a(\lambda x)$ for every $\lambda
\in \mathbb C, a\in \mathcal A$ and $x\in {\mathcal H}$. If there
exists a mapping $\langle \cdot , \cdot \rangle : {\mathcal
H}\times{\mathcal H} \rightarrow \mathcal A$ with the properties

(1) $\langle x, x\rangle \geq 0$ for every $x\in {\mathcal H}$,

(2) $\langle x, x\rangle =0$ if and only if $x=0$,

(3) $\langle x, y\rangle =\langle y, x\rangle ^*$ for every $x,
y\in {\mathcal H}$,

(4) $\langle ax, y\rangle =a\langle x, y\rangle $ for every $a\in
\mathcal A$, every $x, y\in {\mathcal H}$,

(5) $\langle x+y, z\rangle =\langle x, z\rangle +\langle y,
z\rangle $ for every $x, y, z\in {\mathcal H}$.

Then the pair $\{ {\mathcal H}, \langle \cdot , \cdot \rangle \} $
is called a (left)- pre-Hilbert $\mathcal A$-module. The map
$\langle \cdot , \cdot \rangle $ is said to be an $\mathcal
A$-valued inner product. If the pre-Hilbert $\mathcal A$-module
$\{ {\mathcal H}, \langle \cdot , \cdot \rangle \} $ is complete
with respect to the norm $\| x\| =\| \langle x, x\rangle \| ^\frac
{1}{2}$ then it is called a \textit {Hilbert $\mathcal A$-module}.

A Hilbert $\mathcal A$-module ${\mathcal H}$ is (algebraically)
\textit {finitely generated} if there exists a finite set $\{ x_1,
\dots , x_n\} \subseteq {\mathcal H}$ such that every element
$x\in {\mathcal H}$ can be expressed as an $\mathcal A$-linear
combination $x=\sum _{i=1}^{n}a_i, a_i\in \mathcal A$. A Hilbert
$\mathcal A$-module is \textit {countably generated} if there
exists a countable set of generator.

It should be mentioned that by no means all results of Hilbert
space theory can be simply generalized to the situation of Hilbert
$C^*$-modules. Fist of all, the analogue of the Riesz
representation theorem for bounded $\mathcal A$-linear mapping is
not valid for $\mathcal H $. Secondly, the bounded $\mathcal
A$-linear operator on $\mathcal H $ may not have an adjoint
operator. Thirdly, the Hilbert $\mathcal A$-submodule ${\mathcal
I}$ of the Hilbert $\mathcal A$-module $\mathcal H $ is not a
direct summand. Let $\mathcal H$ be a Hilbert $\mathcal A$-module
over a unital $C^*$-algebra $\mathcal A$. The set of all bounded
$\mathcal A$-linear operators on $\mathcal H$ is denoted by
$End_{\mathcal A}({\mathcal H})$, and the set of all adjointable
bounded $\mathcal A$-linear operators on $\mathcal H$ is denoted
by $End_{\mathcal A}^*({\mathcal H})$.

A $C^*$-algebra $\mathcal M$ is called a \textit {$W^*$-algebra}
if it is a dual space as a Banach space, i.e. if there exists a
Banach space ${\mathcal M}_*$ such that $({\mathcal
M}_*)^*=\mathcal M$. We also call ${\mathcal M}_*$ the predual of
$\mathcal M$. It should mention here that $End_{\mathcal
A}^*(l^2({\mathcal A}))$ is not a $W^*$-algebra in general. A
$W^*$-algebra $\mathcal M$ is said to be \textit {finite} if its
identity is finite. Equivalently, $\mathcal M$ is finite if and
only if every isometry in $\mathcal M$ is unitary.

\begin{definition}
Let $\mathcal A$ be a unital $C^*$-algebra and $\mathbb J$ be a
finite or countable index set. A sequence $\{ x_j\}_{j\in {\mathbb
J}}$ of elements in a Hilbert $\mathcal A$-module ${\mathcal H}$
is said to be a (standard) {\it  frame} if there exist two
constants $C, D>0$ such that
\begin{eqnarray}
C\cdot \langle x, x\rangle \leq \sum _{j\in {\mathbb J}}\langle x,
x_j\rangle \langle x_j, x\rangle \leq D\cdot \langle x, x\rangle
\label {mofr}
\end{eqnarray}
for every $x\in {\mathcal H}$, where the sum in the middle of the
inequality is convergent in norm. The optimal constants (i.e.
maximal for $C$ and minimal for $D$) are called frame bounds.

The frame $\{ x_j\}_{j\in {\mathbb J}}$ is said to be tight frame
if $C=D$, and said to be Parseval if $C=D=1$.
\end{definition}

Note that not every Hilbert $C^*$-module has an orthonormal basis.
Though any countably generated Hilbert $C^*$-module admits a
frame, there are countably generated Hilbert $C^*$-modules that
contain no orthonormal basis even no orthogonal Riesz basis (see
Example 3.4 in \cite {frla273}).

The main property of frames for Hilbert spaces is the existence of
the reconstruction formula that allows a simple standard
decomposition of every element of the spaces with respect to the
frame. For standard frames we have the following reconstruction
formula.

\begin{theorem}
(\cite{frla273}) Let $\{ x_j\}_{j\in {\mathbb J}}$ be a standard
frame in a finitely or countably generated Hilbert $A$-module
$\mathcal H$ over a unital $C^*$-algebra $\mathcal A$. Then there
exists a unique operator $S\in End_{\mathcal A}^*(\mathcal H)$
such that
\begin{eqnarray*}
x=\sum _{j\in {\mathbb J}}\langle x, S(x_j)\rangle x_j
\end{eqnarray*}
for every $x\in \mathcal H$. The operator can be explicitly given
by the formula $S=T^*T$ for any adjointable invertible bounded
operator $T$ mapping $\mathcal H$ onto some other Hilbert
$\mathcal A$-module ${\mathcal K}$ and realizing $\{ T(x_j): {j\in
{\mathbb J}} \}$ to be a standard Parseval frame in $\mathcal K$.
\end{theorem}

Let ${\mathcal S}\subseteq End_{\mathcal A}({\mathcal H})$, we
denote its {\it commutant} $\{ A\in End_{\mathcal A}({\mathcal
H}): AS=SA, S\in {\mathcal S}\} $ by ${\mathcal S}^{\prime }$. For
a non-empty set $\mathcal U$ and a unital $C^{*}$-algebra
$\mathcal A$,  let  $l_{\mathcal U}^2({\mathcal A})$ be the
Hilbert $\mathcal A$-module defined by
$$l_{\mathcal U}^2({\mathcal A})=\{ \{ a_U\} _{U\in
{\mathcal U}}\subseteq {\mathcal A}: \sum _{U\in \mathcal
U}a_Ua_U^* \ \textrm {convergences in } \| \cdot \| \} .$$ Let $\{
\chi _U \} _{U\in\mathcal U}$ denote the standard orthonormal
basis of $l_{\mathcal U}^2({\mathcal A})$, where $\chi _U$ takes
value $1_{\mathcal A}$ at $U$ and $0_{\mathcal A}$ at everywhere
else. In the case when $\mathcal U$ is a group, we define for each
$U\in \mathcal U$,
$$L_U\chi _V=\chi _{UV} \ \ \textrm {and} \ \ R_U\chi _V=\chi _{VU^{-1}}.$$
Note that $L_U^{-1}=L_U^*=L_{U^*}$ and $R_U^{-1}=R_U^*=R_{U^*}$.
Here $L$ and $R$ are the {\it left} and {\it right regular
representations} of $\mathcal U$.

Let $\mathcal H$ be a Hilbert $\mathcal A$-module over a unital
$C^*$-algebra $\mathcal A$. A vector $\psi $ in $\mathcal H$ is
called a \textit {wandering vector} for a unitary group $\mathcal
U$ on $\mathcal H$  if ${\mathcal U}\psi =\{U\psi:  U\in {\mathcal
U}\} $ is an orthonormal set. If ${\mathcal U}\psi $ is an
orthonormal basis for $\mathcal H$, then $\psi $ is called a
\textit {complete wandering vector} for $\mathcal U$. Similarly, a
vector $\psi $ is called a \textit {Parseval frame vector} (resp.
\textit {frame vector with bounds $C$ and $D$}, or \textit {Bessel
sequence vector with bound $D$}) for a unitary group $\mathcal U$
if ${\mathcal U}x$  forms a Parseval frame (resp. frame with
bounds $C$ and $D$, or Bessel sequence with bound $D$) for
$\overline{span}({\mathcal U}x)$.  Moreover, $x$ is called a
\textit {complete Parseval frame vector} (resp. \textit {complete
frame vector with bounds $C$ and $D$}, or \textit {complete Bessel
sequence with bound $D$}) when ${\mathcal U}x$ is a Parseval frame
(resp. frame with bounds $C$ and $D$, or Bessel sequence with
bound $D$) for $\mathcal H$.

The following simple lemma will be used in the proof of Theorem
\ref{the-3.4}.

\begin{lemma}\label{dilation}
Let $\mathcal G$ be a unitary group on a finitely or countably
generated Hilbert $\mathcal A$-module $\mathcal H$ over a unital
$C^*$-algebra $\mathcal A$. If $\mathcal G$ admits a complete
Parseval frame vector $\eta$, then $\mathcal G$ is unitarily
equivalent to $\{L_{U}|_{\mathcal K}: U\in\mathcal{U}\}$, where
${\mathcal K}= T(\mathcal H)$ and $T: {\mathcal H}\rightarrow
l_{\mathcal G}^2(\mathcal A)$ be the analysis operator defined by
$T(x)=\sum _{U\in \mathcal G}\langle  x, U\eta \rangle \chi _U.$
\end{lemma}

\begin{proof} It is easy to check that  $T$ is an adjointable isometry.
By Theorem 15.3.5 and Theorem 15.3.8 in \cite{weol1993}, we have
that
$$\l_{\mathcal G}^2(\mathcal A)=(T({\mathcal H}))^{\perp }\oplus T({\mathcal H}).$$
Hence we have the orthogonal projection $P$ from $l_{\mathcal
G}^2(\mathcal A)$ onto $T({\mathcal H})$.

For each $V\in \mathcal G$, we have
\begin{eqnarray*}
L_UT(V\eta )&=&L_U(\sum_{W\in \mathcal G}\langle  V\eta , W\eta \rangle   \chi _W)=\sum_{W\in \mathcal G}\langle  V\eta , W\eta \rangle   \chi _{UW}\\
&=&\sum_{W\in \mathcal G}\langle  UV\eta , UW\eta \rangle   \chi _{UW}=\sum_{W\in \mathcal G}\langle  UV\eta , W\eta \rangle   \chi _W\\
&=&TU(V\eta ).
\end{eqnarray*}
Thus $L_UT=TU$.
\end{proof}

We remark that the orthogonal projection from $l_{\mathcal
G}^2(\mathcal A)$ onto $T(\mathcal H)$ is in the commutant of
$\{L_{U}: U\in\mathcal G\}$ and satisfies $T\eta = P\chi_{I}$.
This also implies the so-called {\it dilation property} meaning
that there exists a Hilbert $\mathcal A$-module $\tilde {\mathcal
H}\supseteq \mathcal H$ and a unitary group $\tilde {\mathcal G}$
on $\tilde {\mathcal H}$ such that $\tilde {\mathcal G}$ has
complete wandering vectors in $\tilde {\mathcal H}$, $\mathcal H$
is an invariant subspace of $\tilde {\mathcal G}$ such that
$\tilde {\mathcal G}|_{\mathcal H}=\mathcal G$, and the map
$G\mapsto G|_{\mathcal H}$ is a group isomorphism from $\tilde
{\mathcal G}$ onto $\mathcal G$.

\section {Frame Vector Parameterizations}

In \cite {dala1998} the set of all  wandering vectors for a
unitary group was parameterized by the set of unitary operators in
its commutant. However,  unlike the wandering vector case, it was
shown in \cite {hala2000} that the set of all the Parseval frame
vectors for a unitary group can not be parameterized by the set of
all the unitary operators in the communtant of the unitary group.
This means that the Parseval frame vectors for a representation of
a countable group are not necessarily unitarily equivalent.
However, this set can be parameterized by the set of all the
unitary operators in the von Neumann algebra generated by the
representation (\cite{ha3329, hala2000}). This turns out to be a
very useful result in Gabor analysis (cf. \cite{ha3329,
hagala2000}). Although it remains a question whether this result
is still valid in the Hilbert $C^{*}$-module setting, we will
prove in this section that this result holds when the underlying
$C^{*}$-algebra is a commutative $W^*$-algebra. Even in this
commutative case, a lot more extra work and care are needed in
order to prove this generalization.

\begin{lemma}\label{lem-3.3}
Let $\mathcal G$ be a unitary group on a Hilbert ${\mathcal
A}$-module over a commutative unital $C^*$-algebra $\mathcal A$,
then
$${\mathcal M}={\mathcal N}^{\prime }=\{ R_U: U\in {\mathcal G}\} ^{\prime} \ \
{\textrm and }\ \ {\mathcal N}={\mathcal M}^{\prime }=\{ L_U: U\in {\mathcal G}\} ^{\prime },$$
where ${\mathcal M}=\{ L_U: U\in {\mathcal G}\} ^{\prime \prime}$
and ${\mathcal N}=\{ R_U: U\in {\mathcal G}\} ^{\prime \prime}$.
\end{lemma}
\begin{proof}
Note that $R_UL_V=L_VR_U$ holds for any $U, V\in \mathcal G$.
Therefore to prove this lemma it suffices to show that $TS=ST$ for
arbitrary $T\in {\mathcal M}^{\prime }$ and $S\in {\mathcal
N}^{\prime }$.

Suppose that $$T\chi _I=\sum _{U\in {\mathcal G}}a_U\chi _U\ \
\textrm {and}\ \ S\chi _I=\sum _{U\in {\mathcal G}}b_U\chi _U$$
for some $a_U, b_U\in \mathcal A$.

Now for any $V\in \mathcal G$, on one hand, we have
\begin{eqnarray*}
ST\chi _V&=&STL_V\chi _I=SL_VT\chi _I \\
&=&SL_V(\sum _{U\in {\mathcal G}}a_U\chi _U)=S(\sum _{U\in {\mathcal G}}a_U\chi _{VU})\\
&=&S(\sum _{U\in {\mathcal G}}a_UR_{(VU)^{-1}}\chi _{I})=\sum _{U\in {\mathcal G}}a_UR_{(VU)^{-1}}S\chi _I\\
&=&\sum _{U\in {\mathcal G}}a_UR_{(VU)^{-1}}(\sum _{W\in {\mathcal
G}}b_W\chi _W)=\sum _{U, W\in {\mathcal G}}a_Ub_W\chi _{WVU}.
\end{eqnarray*}
On the other hand
\begin{eqnarray*}
TS\chi _V&=&TSR_{V^{-1}}\chi _I=TR_{V^{-1}}S\chi _I\\
&=&TR_{V^{-1}}(\sum _{W\in {\mathcal G}}b_W\chi _W)=T(\sum _{W\in {\mathcal G}}b_W\chi _{WV})\\
&=&T(\sum _{W\in {\mathcal G}}b_WL_{WV}\chi _{I})=\sum _{W\in {\mathcal G}}b_WL_{WV}T\chi _I\\
&=&\sum _{W\in {\mathcal G}}b_WL_{WV}(\sum _{U\in {\mathcal
G}}a_U\chi _U)=\sum _{U, W\in {\mathcal G}}b_Wa_U\chi _{WVU}.
\end{eqnarray*}
Since $\mathcal A$ is commutative, it follows that $ST\chi
_V=TS\chi _V$, and so $ST=TS$.
\end{proof}

We now define a natural conjugate ${\mathcal A}$-linear
isomorphism $\pi $ from $\mathcal M$ onto $ {\mathcal M}^{\prime
}=\mathcal N$ by
$$\pi (A)B\chi _I=BA^*\chi _I, \ \  \forall A, B\in {\mathcal M}.$$
In particular, $\pi (A)\chi _I=A^*\chi _I$.

Now we are in a position to prove the parameterization of complete
Parseval frame vectors for unitary groups.

\begin{theorem}\label{the-3.4}
Let $\mathcal G$ be a unitary group on a finitely or countably
generated Hilbert $\mathcal A$-module $\mathcal H$ over a
commutative  $W^*$-algebra $\mathcal A$ such that $l_{\mathcal
G}^2({\mathcal A})$ is self-dual. Suppose that $\eta\in \mathcal
H$ be a complete Parseval frame vector for $\mathcal G$. For $\xi
$ in $\mathcal H$ we have

(1) $\xi $ is a complete Parseval frame vector for $\mathcal G$ if
and only if there exists a unitary operator $A\in {\mathcal
G}^{\prime \prime }$ such that $\xi =A\eta $.

(2) $\xi $ is a complete   frame vector for $\mathcal G$ if and
only if there exists an invertible and adjointable  operator $A\in
{\mathcal G}^{\prime \prime }$ such that $\xi =A\eta $.

(3) $\xi $ is a complete  Bessel sequence vector for $\mathcal G$
if and only if there exists an adjointable  operator $A\in
{\mathcal G}^{\prime \prime }$ such that $\xi =A\eta $.
\end{theorem}
\begin{proof} We will prove $(1)$. The proof of $(2)$ and $(3)$ is similar and we leave
it to the interested readers.

By Lemma \ref{dilation}, we can assume that ${\mathcal G}=\{
L_U|_{Rang(P)}, U\in {\mathcal G}\} $ and $\eta =P\chi _I$, where
$P$ is an orthogonal projection in the commutant of $\{L_{U}:
U\in\mathcal G\}$

Let ${\mathcal M}=\{ L_U: U\in {\mathcal G}\}''$

First assume that there exists a unitary operator $A\in {\mathcal
G} ^{\prime \prime }$ such that $\xi =A\eta$.

We now show that $A\eta$ is a complete Parseval frame vector for
$\mathcal G$. For any $x\in Rang(P)$, we have
\begin{eqnarray*}
& &\sum _{U\in \mathcal {G}}\langle   x, UA\eta \rangle   \langle   UA\eta, x\rangle
=\sum _{U\in \mathcal G}\langle   x, L_UPA\eta\rangle   \langle   L_UPA\eta , x\rangle   \\
&=&\sum _{U\in \mathcal G}\langle   x, L_UPAP\chi _I\rangle   \langle   L_UPAP\chi _I, x\rangle
=\sum _{U\in \mathcal G}\langle   x, L_UPA\chi _I\rangle   \langle   L_UPA\chi _I, x\rangle   \\
&=&\sum _{U\in \mathcal G}\langle   x, PL_UA\chi _I\rangle   \langle   PL_UA\chi _I, x\rangle
=\sum _{U\in \mathcal G}\langle   Px, L_UA\chi _I\rangle   \langle   L_UA\chi _I, Px\rangle   \\
&=&\sum _{U\in \mathcal G}\langle   x, L_U\pi (A^*)\chi _I\rangle   \langle   L_U\pi (A^*)\chi _I, x\rangle   \\
&=&\sum _{U\in \mathcal G}\langle   x, \pi (A^*)L_U\chi _I\rangle   \langle   \pi (A^*)L_U\chi _I, x\rangle   \\
&=&\sum _{U\in \mathcal G}\langle   (\pi (A^*))^*x,  \chi _U\rangle   \langle   \chi _U, (\pi (A^*))^*x\rangle   \\
&=&\langle   (\pi (A^*))^*x, (\pi (A^*))^*x\rangle  =\langle   x,
x\rangle   ,
\end{eqnarray*}
where in the seventh equality we use that fact $\pi
(A^*)L_U=L_U\pi (A^*)$, and in the last equality we use that fact
that $\pi (A^*)$ is unitary. Therefore  $A\eta $ is a complete
Parseval frame vector for $\mathcal G$.

Now let $\xi \in Rang (P)$ be a complete Parseval frame vector for
$\mathcal G$. We want to find a unitary operator $A\in {\mathcal
G} ^{\prime \prime }$ such that $\xi =A\eta$.

To this aim, we first define an operator $B: l_{\mathcal
G}^2({\mathcal A})\rightarrow l_{\mathcal G}^2({\mathcal A})$ by
$$\chi _U\longmapsto L_U\xi,  \ \ U\in {\mathcal G}.$$

One can check that $B$ is an adjointable operator and $B^*\chi
_V=\sum _{W\in {\mathcal G}}\langle   L_{W^{-1}}L_V\eta,
\xi\rangle   \chi _W$ for any $V\in \mathcal G$.

Now for any $U, V\in \mathcal G$, we see that
\begin{eqnarray*}
& &\langle   (BB^*-P)\chi _U, \chi _V\rangle \\
&=&\langle   \sum _{W\in \mathcal G}\langle   L_{W^{-1}}L_U\eta,
\xi \rangle   \chi _W, \sum _{S\in \mathcal G}\langle
L_{S^{-1}}L_V\eta , \xi \rangle
 \chi _S\rangle-\langle   T_{\eta }U\eta , \chi _V\rangle   \\
&=&\sum _{W\in \mathcal G}\langle   L_{W^{-1}}L_U\eta , \xi \rangle   \langle   \xi , L_{W^{-1}}L_V\eta \rangle
-\langle   \sum _{W\in \mathcal G}\langle   U\eta , W\eta \rangle   \chi _W, \chi _V\rangle   \\
&=&\sum _{W\in \mathcal G}\langle   L_U\eta , L_W\xi \rangle   \langle   L_W\xi , L_V\eta \rangle   -\langle   U\eta , V\eta \rangle   \\
&=&\langle   L_U\eta , L_V\eta \rangle   -\langle   U\eta , V\eta \rangle   =\langle
L_UP\chi _I, L_VP\chi _I\rangle   -\langle   U\eta , V\eta \rangle   \\
&=&\langle   L_UT_{\eta }\eta , L_VT_{\eta }\eta \rangle
-\langle   U\eta , V\eta \rangle   =\langle   T_{\eta }U\eta ,
T_{\eta }V\eta \rangle   -\langle   U\eta , V\eta \rangle =0,
\end{eqnarray*}
this leads to the fact that $P=BB^*$.

From $$BL_U\chi _V=B\chi _{UV}=L_{UV}\xi =L_UL_V \xi =L_UB\chi
_V,$$ we see that $B\in {\mathcal M}^{\prime }$. Hence $B$ is a
partial isometry in ${\mathcal M}^{\prime }$.

Let $Q=B^*B$, then $P$ and $Q$ are equivalent projections in
${\mathcal M}^{\prime }$.

Since $l_{\mathcal G}^2({\mathcal A})$ is self-dual, by \cite
{pa443}, $End_{\mathcal A}^*(l_{\mathcal G}^2({\mathcal A}))$ is a
$W^*$-algebra. Let $(End_{\mathcal A}^*(l_{\mathcal G}^2({\mathcal
A})))_*$ be its predual. One can check that $\mathcal M$ and
${\mathcal M}^{\prime }$ are  $\sigma (End_{\mathcal
A}^*(l_{\mathcal G}^2({\mathcal A})), (End_{\mathcal
A}^*(l_{\mathcal G}^2({\mathcal A})))_*)$-closed in $End_{\mathcal
A}^*(l_{\mathcal G}^2({\mathcal A}))$, and so both $\mathcal M$
and ${\mathcal M}^{\prime }$ are $W^*$-algebras (see
\cite{sa1998}).

\vskip 2mm {\bf Claim}. {\rm $\mathcal M$ and ${\mathcal
M}^{\prime }$ are finite $W^*$-algebras.} \vskip 2mm We now define
$\phi : {\mathcal M}\rightarrow {\mathcal A}$ by
$$\phi (A)=\langle   A\chi _I, \chi _I\rangle   , \ \  \forall A\in {\mathcal M}.$$

We want to show that $\phi $ is a faithful $\mathcal A$-valued
trace for $\mathcal M$.

Since $\overline {span \{ L_U\chi _I, U\in {\mathcal
G}\}}=l_{\mathcal G}^2({\mathcal A})$, for any $A, B\in \mathcal
M$, we have
$$A\chi _I=\lim _{n}A_n\chi _I \ \  {\rm and }\ \ B\chi _I=\lim _{n}B_n\chi _I ,$$
where
$$A_n\chi _I=\sum _{i=1}^{k_n}a_i^{(n)}L_{{V_i}^{(n)}}\chi _I \ \ {\rm and }\ \ B_n\chi _I=\sum _{j=1}^{l_n}b_j^{(n)}L_{{W_j}^{(n)}}\chi _I $$
for some $a_i^{(n)}, b_j^{(n)}\in \mathcal A$ and ${V_i}^{(n)},
{W_j}^{(n)}\in \mathcal G$.

Then
$$\phi (AB)=\langle   AB\chi _I, \chi _I\rangle   =\lim _{m}\lim _{n}\langle   \sum _{j=1}^{l_m}\sum _{i=1}^{k_n}b_j^{(m)}a_i^{(n)}L_{{W_j}^{(m)}}L_{{V_i}^{(n)}}\chi _I, \chi _I\rangle   .$$
While
$$ \phi (BA)=\lim _{n}\lim _{m}\langle   \sum _{i=1}^{k_n}\sum _{j=1}^{l_m}a_i^{(n)}b_j^{(m)}L_{{V_i}^{(n)}}L_{{W_j}^{(m)}}\chi _I, \chi _I\rangle   .$$

Note that
$$\langle   L_{{W_j}^{(m)}}L_{{V_i}^{(n)}}\chi _I, \chi _I\rangle   =\langle   L_{{V_i}^{(n)}}L_{{W_j}^{(m)}}\chi _I, \chi _I\rangle   .$$

Therefore $\phi (AB)=\phi (BA)$.

If $A\in \mathcal M$ is positive and $\phi (A)=0$, then
$$\langle   A^{\frac {1}{2}}\chi _I, A^{\frac {1}{2}}\chi _I\rangle   =\langle   A\chi _I, \chi _I\rangle =\phi (A)=0.$$
Thus $A^{\frac {1}{2}}\chi _I=0$.

Now for any $U\in \mathcal G$, we have
$$A^{\frac {1}{2}}\chi _U=A^{\frac {1}{2}}R_U\chi _I=R_UA^{\frac {1}{2}}\chi _I=0.$$
Therefore $A^{\frac {1}{2}}=0$, and so $A=0$. Similarly, by using
Lemma \ref{lem-3.3} we can prove that ${\mathcal M}^{\prime }$ is
also finite.

It follows from Proposition 2.4.2 in \cite {sa1998} that $I-P$ and
$I-Q$ are equivalent projections in ${\mathcal M}^{\prime }$.
Therefore there exists a partial isometry $C\in {\mathcal
M}^{\prime }$ such that $CC^*=I-P$ and $C^*C=I-Q$.

Let $T=B+C$. Then $T$ is a unitary operator in ${\mathcal
M}^{\prime }$, and so $A=(\pi ^{-1}(T))^*$ is a unitary operator
in $\mathcal M$.

In order to  complete the proof it remains to prove that $A\tilde
{\eta }=\tilde {\xi }$.

In fact,
\begin{eqnarray*}
A\eta &=&(\pi ^{-1}(T))^*P\chi _I=P(\pi ^{-1}(T))^*\chi _I\\
&=&P\pi (\pi ^{-1}(T))\chi _I=PT\chi _I\\
&=&P(B+C)\chi _I=PB\chi _I+PC\chi _I\\
&=&P\xi =\xi ,
\end{eqnarray*}
which completes the proof.
\end{proof}

The following result follows immediately from Theorem
\ref{the-3.4} and the fact the set of all the unitary elements in
any $W^{*}$-algebra is path-connected in norm.

\begin{corollary}
Let $\mathcal G$ be a unitary group on a finitely or countably
generated Hilbert $\mathcal A$-module $\mathcal H$ over a
commutative $W^*$-algebra $\mathcal A$ such that $l_{\mathcal
G}^2({\mathcal A})$ is self-dual, then the set of all Paserval
frame vectors for $\mathcal G$ is path-connected.

\end{corollary}

\section {Parseval Frame Approximations}

In the Hilbert space frame setting, the original work on symmetric
orthogonalization was done by L$\ddot{\textrm{o}}$widin \cite
{lo185} in the late 1970's. The concept of symmetric approximation
of frames by Parseval frame was introduced in \cite {frpati777} to
extend the symmetric orthogonalization of bases by orthogonal
bases in Hilbert spaces. The existence and the uniqueness results
for the symmetric approximation of frames by Parseval frames were
obtained in \cite{frpati777}. Following their definition, a
Parseval frame $\{ y_j\} _{j=1}^{\infty }$ is said to be a \textit
{symmetric approximation} of frame $\{ x_j\} _{j=1}^{\infty }$ in
Hilbert space $H$ if it is similar to $\{ x_j\} _{j=1}^{\infty }$
and
\begin{eqnarray}
\sum _{j=1}^{\infty}\| z _j-x_j\| ^2\geq \sum
_{j=1}^{\infty}\|y_j-x_j\| ^2 \label {app}
\end{eqnarray}
is valid for all Parseval frames $\{ z_j\} _{j=1}^{\infty }$ of
$H$ that are similar to $\{ x_j\} _{j=1}^{\infty }$.

Observed by the first author (\cite{ha3329}, \cite{ ha78}) that in
some situations the symmetric approximation fails to work when the
underlying Hilbert space is infinite dimensional since if we
restrict ourselves to the frames induced by a unitary system then
the summation in (\ref{app}) is always infinite when the given
frame is not Parseval.  Instead of using the symmetric
approximations to consider the frames generated by a collection of
unitary transformations and some window functions, it was proposed
to approximate the frame generator by Parseval frame generators.
Existence and uniqueness results for such a best approximation
were obtained in \cite{ha3329, ha78}. We will extend this result
to Hilbert $C^*$-module frames when the underlying $C^{*}$-algebra
is commutative. It remains open whether this is true when the
underlying $C^{*}$-algebra is non-commutative.

Following (\cite {ha78}) we first give the following definition.

\begin{definition}
Let $\Phi =( \phi _1, \dots , \phi _N)$ be a multi-frame generator
for a unitary system $\mathcal U$. Then a Parseval multi-frame
generator $\Psi =(\psi _1, \dots , \psi _N)$ for $\mathcal U$ is
called a best Parseval multi-frame approximation for $\Phi $ if
the inequality
$$\sum _{k=1}^{N}\langle \phi _k-\psi _k, \phi _k-\psi _k\rangle \leq \sum _{k=1}^{N}\langle \phi _k-\xi _k, \phi _k-\xi _k\rangle $$
is valid for all the Parseval multi-frame generator $\Xi =(\xi _1,
\dots , \xi _N)$ for $\mathcal U$.
\end{definition}

Let $\Phi =\{ \phi _1, \phi _2, \dots , \phi _N\} $ be a
multi-frame generator for a unitary system $\mathcal U$ on a
finitely or countably generated Hilbert $\mathcal A$-module
$\mathcal H$ over a unital $C^*$-algebra $\mathcal A$. We use
$T_{\Phi }$ to denote the analysis operator from $\mathcal H$ to
$l_{{\mathcal U}\times \{ 1, 2, \dots , N\}}^2({\mathcal A})$
defined by
$$T_{\Phi }x=\sum _{j=1}^{N}\sum _{U\in \mathcal U}\langle  x, U\phi _j\rangle   \chi _{(U, j)}, \ \  \forall x\in \mathcal H,$$
where $\{ \chi _{(U, j)}: U\in {\mathcal U}, j=1, 2, \dots, N\} $
is the standard orthonormal basis for $l_{{\mathcal U}\times \{ 1,
2, \dots , N\}}^2({\mathcal A})$.

Note that $T_{\Phi }$ is adjointable and its adjoint operator
satisfying
$$T_{\Phi }^*\chi _{(U, j)}=U\phi _j, \ \  U\in {\mathcal U},\ j=1, 2, \dots ,N.$$

\begin{lemma}\label{magic}
Let $\mathcal G$ be a unitary group on a Hilbert $\mathcal
A$-module $\mathcal H$ over a commutative $C^*$-algebra $\mathcal
A$. Suppose that $\Phi =\{ \phi _1, \phi _2, \dots , \phi _N\} $
and $\Psi =\{ \psi _1, \psi _2, \dots , \psi _N\} $ be two
multi-frame generators for $\mathcal G$, then
$$\sum _{k=1}^{N}\langle  \phi _k, \phi _k\rangle  =\sum _{k=1}^{N}\langle  \psi _k, \psi _k\rangle  .$$
\end{lemma}

\begin{proof} We compute
\begin{eqnarray*}
\sum _{k=1}^{N}\langle  \phi _k, \phi _k\rangle  &=&\sum _{k=1}^{N}\sum _{j=1}^{N}\sum _{U\in \mathcal G}\langle  \phi _k, U\psi _j\rangle  \langle  U\psi _j, \phi _k\rangle  \\
&=&\sum _{k=1}^{N}\sum _{j=1}^{N}\sum _{U\in \mathcal G}\langle  U^*\phi _k, \psi _j\rangle  \langle  \psi _j, U^*\phi _k\rangle  \\
&=&\sum _{j=1}^{N}\sum _{k=1}^{N}\sum _{U\in \mathcal G}\langle  \psi _j, U^*\phi _k\rangle  \langle  U^*\phi _k, \psi _j\rangle  \\
&=&\sum _{j=1}^{N}\langle  \psi _j, \psi _j\rangle  .
\end{eqnarray*}
\end{proof}

\begin{theorem} \label{mul-app}Let $\mathcal G$ be a unitary group on a finitely or countably generated Hilbert $\mathcal A$-module
$\mathcal H$ over a commutative unital $C^*$-algebra $\mathcal A$.
Suppose that $\Phi =\{ \phi _1, \phi _2, \dots , \phi _N\} $ is a
multi-frame generator for $\mathcal G$. Then $S^{\frac {1}{2}}\Phi
$ is the unique best Parseval multi-frame approximation for $\Phi
$, where $S$ is the frame operator for the multi-frame $\{ U\phi
_j: j=1, \dots , N, U\in \mathcal G \}$.
\end{theorem}

\begin{proof} We first show that $S\in {\mathcal G}^{\prime }$.

For arbitrary $V\in \mathcal G$ and $x\in \mathcal H$ we have
\begin{eqnarray*}
SVx&=&\sum _{k=1}^{N}\sum _{U\in \mathcal G}\langle  Vx, U\phi _k\rangle   U\phi _k\\
&=&\sum _{k=1}^{N}\sum _{U\in \mathcal G}\langle  x, V^*U\phi  _k\rangle   U\phi _k\\
&=&V(\sum _{k=1}^{N}\sum _{U\in \mathcal G}\langle  x, V^*U\phi _k\rangle   V^*U\phi _k)\\
&=&V(\sum _{k=1}^{N}\sum _{U\in \mathcal G}\langle  x, U\phi _k\rangle   U\phi _k)\\
&=&VSx.
\end{eqnarray*}
This shows that $S\in {\mathcal G}^{\prime }$.

Since $End_{\mathcal A}^*(\mathcal H)$ is a $C^*$-algebra, by the
spectral decomposition for positive elements in $C^*$-algebra, we
can infer that $S^{-\frac {1}{2}}, S^{-\frac {1}{4}}\in {\mathcal
G}^{\prime}$. Therefore $\{ S^{-\frac {1}{2}}\phi _1,  S^{-\frac
{1}{2}}\phi _2,  \dots , S^{-\frac {1}{2}}\phi _N\} $ is a
complete Parseval multi-frame generator for $\mathcal G$.

Let $\Psi =\{ \psi _1, \psi _2, \dots ,\psi _N\} $ be any Parseval
multi-frame generator for $\mathcal G$. We claim that
$$\sum _{k=1}^{N}\langle  T_{S^{-\frac {1}{2}}\Phi}S^{-\frac {1}{4}}\phi _k, T_{\Psi }S^{-\frac {1}{4}}\phi _k\rangle   =\sum _{k=1}^{N}\langle \psi _k, \phi _k\rangle   ,$$
where $T_{S^{-\frac {1}{2}}\Phi }$ and $T_{\Psi }$ are the
analysis operators with respect to the Parseval multi-frame
generators $S^{-\frac {1}{2}}\Phi $ and $\Psi $ respectively.

We compute
\begin{eqnarray*}
& &\sum _{k=1}^{N}\langle  T_{S^{-\frac {1}{2}}\Phi }S^{-\frac {1}{4}}\phi _k, T_{\Psi }S^{-\frac {1}{4}}\phi _k\rangle   \\
&=&\sum _{k=1}^{N}\langle \sum _{j=1}^{N}\sum _{U\in \mathcal G}\langle  S^{-\frac {1}{4}}\phi _k, US^{-\frac {1}{2}}\phi _j\rangle   \chi _{(U, j)}, \sum _{i=1}^{N}\sum _{V\in \mathcal G}\langle  S^{-\frac {1}{4}}\phi _k, V\psi _i\rangle   \chi _{(V, i)}\rangle   \\
&=&\sum _{k=1}^{N}\sum _{j=1}^{N}\sum _{U\in \mathcal G}\langle  S^{-\frac {1}{4}}\phi _k, US^{-\frac {1}{2}}\phi _j\rangle   \langle  U\psi _j, S^{-\frac {1}{4}}\phi _k\rangle   \\
&=&\sum _{k=1}^{N}\sum _{j=1}^{N}\sum _{U\in \mathcal G}\langle  U\psi _j, S^{-\frac {1}{4}}\phi _k\rangle   \langle  S^{-\frac {1}{4}}\phi _k, US^{-\frac {1}{2}}\phi _j\rangle   \\
&=&\sum _{j=1}^{N}\sum _{k=1}^{N}\sum _{U\in \mathcal G}\langle  S^{\frac {1}{4}}\psi _j, U^*S^{-\frac {1}{2}}\phi _k\rangle   \langle  U^*S^{-\frac {1}{2}}\phi _k, S^{-\frac {1}{4}}\phi _j\rangle   \\
&=&\sum _{j=1}^{N}\langle  S^{\frac {1}{4}}\psi _j, S^{-\frac
{1}{4}}\phi _j\rangle  =\sum _{j=1}^{N}\langle \psi _j, \phi
_j\rangle .
\end{eqnarray*}

We now prove that $S^{\frac {1}{2}}\Phi $ is a best Parseval
multi--frame approximation for $\Phi $. We need to show that
$$\sum _{k=1}^{N}\langle \psi _k-\phi _k, \psi _k-\phi _k\rangle   \geq \sum _{k=1}^{N}\langle  S^{-\frac {1}{2}}\phi _k-\phi _k, S^{-\frac {1}{2}}\phi _k-\phi _k\rangle   .$$

By Lemma \ref{magic}, it suffices to prove that
$$\sum _{k=1}^{N}(\langle  S^{-\frac {1}{2}}\phi _k, \phi _k\rangle   +\langle \phi _k, S^{-\frac {1}{2}}\phi _k)-\langle \psi _k, \phi _k\rangle   -\langle \phi _k, \psi _k\rangle   )\geq 0.$$

In fact, we have
\begin{eqnarray*}
& &\sum _{k=1}^{N}(\langle  S^{-\frac {1}{2}}\phi _k, \phi
_k\rangle   +\langle \phi _k, S^{-\frac {1}{2}}\phi _k\rangle
-\langle \psi _k, \phi _k\rangle   -\langle \phi _k, \psi
_k\rangle   )\\ \nonumber
&=&\sum _{k=1}^{N}(\langle  S^{-\frac {1}{4}}\phi _k, S^{-\frac {1}{4}}\phi _k\rangle   +\langle S^{-\frac {1}{4}}\phi _k, S^{-\frac {1}{4}}\phi _k)\rangle   \\
& &-\langle  T_{S^{-\frac {1}{2}}\Phi }S^{-\frac {1}{4}}\phi _k, T_{\Psi }S^{-\frac {1}{4}}\phi _k\rangle   -\langle  T_{\Psi }S^{-\frac {1}{4}}\phi _k, T_{S^{-\frac {1}{2}}\Phi }S^{-\frac {1}{4}}\phi _k\rangle   )\\
&=&\sum _{k=1}^{N}(\langle  T_{S^{-\frac {1}{2}}\Phi }S^{-\frac {1}{4}}\phi _k, T_{S^{-\frac {1}{2}}\Phi }S^{-\frac {1}{4}}\phi _k\rangle   +\langle  T_{\Psi }S^{-\frac {1}{4}}\phi _k, T_{\Psi }S^{-\frac {1}{4}}\phi _k\rangle   \\
& &-\langle  T_{S^{-\frac {1}{2}}\Phi }S^{-\frac {1}{4}}\phi _k, T_{\Psi }S^{-\frac {1}{4}}\phi _k\rangle   -\langle  T_{\Psi }S^{-\frac {1}{4}}\phi _k, T_{S^{-\frac {1}{2}}\Phi }S^{-\frac {1}{4}}\phi _k\rangle   )\\
&=&\sum _{k=1}^{N}\langle (T_{S^{-\frac {1}{2}}\Phi }-T_{\Psi
})S^{-\frac {1}{4}}\phi _k, (T_{S^{-\frac {1}{2}}\Phi }-T_{\Psi
})S^{-\frac {1}{4}}\phi _k\rangle   \geq 0.
\end{eqnarray*}
This implies that $S^{-\frac {1}{2}}\Phi $ is a best Parseval
multi-frame approximation for $\Phi $.

For the uniqueness, assume that $\Xi =\{ \xi _1, \xi _2, \dots ,
\xi _N\}$ be another best Parseval multi-frame approximation for
$\Phi $. Then we have
\begin{equation}
\sum _{k=1}^{N}\langle \xi _k-\phi _k, \xi _k-\phi _k\rangle =\sum
_{k=1}^{N}\langle  S^{-\frac {1}{2}}\phi _k-\phi _k, S^{-\frac
{1}{2}}\phi _k-\phi _k\rangle   .\label{px}
\end{equation}

By Lemma \ref{magic}, we also have
\begin{equation}\sum _{k=1}^{N}\langle \xi _k, \xi _k\rangle   =\sum _{k=1}^{N}\langle  S^{-\frac {1}{2}}\phi _k, S^{-\frac {1}{2}}\phi _k\rangle   .\label{xs}
\end{equation}

Identities (\ref{px}) and (\ref{xs}) imply that
\begin{eqnarray*}
\sum _{k=1}^{N}(\langle \xi _k, \phi _k\rangle   +\langle \phi _k, \xi _k\rangle   )&=&\sum _{k=1}(\langle  S^{-\frac {1}{2}}\phi _k, \phi _k\rangle   +\langle \phi _k, S^{-\frac {1}{2}}\phi _k\rangle   )\\
&=&2\sum _{k=1}^{N}\langle  S^{-\frac {1}{4}}\phi _k, S^{-\frac
{1}{4}}\phi _k\rangle   .
\end{eqnarray*}

We claim that
$$\sum _{k=1}^{N}\langle  S^{\frac {1}{4}}\xi _k, S^{\frac {1}{4}}\xi _k\rangle   =\sum _{k=1}^{N}\langle  S^{-\frac {1}{4}}\phi _k, S^{-\frac {1}{4}}\phi _k\rangle   .$$

In fact,
\begin{eqnarray*}
& &\sum _{k=1}^{N}\langle  S^{\frac {1}{4}}\xi _k, S^{\frac {1}{4}}\xi _k\rangle   \\
&=&\sum _{k=1}^{N}\sum _{j=1}^{N}\sum _{U\in \mathcal G}\langle  S^{\frac {1}{4}}\xi _k, US^{-\frac {1}{2}}\phi _j\rangle   \langle  US^{-\frac {1}{2}}\phi _j, S^{\frac {1}{4}}\xi _k\rangle   \\
&=&\sum _{j=1}^{N}\sum _{k=1}^{N}\sum _{U\in \mathcal G}\langle  S^{-\frac {1}{4}}\phi _j, U^*\xi _k\rangle   \langle  U^*\xi _k, S^{-\frac {1}{4}}\phi _j\rangle   \\
&=&\sum _{j=1}^{N}\langle  S^{-\frac {1}{4}}\phi _j, S^{-\frac
{1}{4}}\phi _j\rangle   .
\end{eqnarray*}
Then we have
 \begin{eqnarray*}
& &\sum _{k=1}^{N}\langle  S^{\frac {1}{4}}\xi _k-S^{-\frac {1}{4}}\phi _k, S^{\frac {1}{4}}\xi _k-S^{-\frac {1}{4}}\phi _k\rangle   \\
&=&\sum _{k=1}^{N}(\langle  S^{\frac {1}{4}}\xi _k, S^{\frac {1}{4}}\xi _k\rangle   -\langle  S^{\frac {1}{4}}\xi _k, S^{-\frac {1}{4}}\phi _k\rangle   \\
& &-\langle  S^{-\frac {1}{4}}\phi _k, S^{\frac {1}{4}}\xi _k\rangle   +\langle  S^{-\frac {1}{4}}\phi _k, S^{-\frac {1}{4}}\phi _k\rangle )\\
&=&\sum _{k=1}^{N}(2\langle  S^{-\frac {1}{4}}\phi _k, S^{-\frac {1}{4}}\phi _k\rangle -\langle \xi _k, \phi _k\rangle -\langle \phi _k, \xi _k\rangle )\\
&=&0.
\end{eqnarray*}
This implies that
$$S^{\frac {1}{4}}\xi _k=S^{-\frac {1}{4}}\phi _k, \ \  k=1, 2, \dots , N.$$
Therefore
$$\xi _k=S^{-\frac {1}{2}}\phi _k, \ \   k=1, 2, \dots , N.$$
i.e. $\Xi=S^{-\frac {1}{2}}\Phi $, as expected.
\end{proof}

\vskip 2mm \textbf {Acknowledgement.} The first author is grateful
to Professor M. Frank for many useful communications and help. The
authors thank the referee for some valuable suggestions in
improving the presentation of the paper.

\bibliographystyle{amsplain}

\end{document}